\documentclass{amsart}
\usepackage{amscd,amssymb}
\usepackage{graphicx}
\numberwithin{equation}{section}

\usepackage{color}

\def\bZ{\Bbb Z}
\def\bC{\Bbb C}

\def\bR{\Bbb R}

\def\tS{\tilde S}

\def\oS{\overline{S}}
\def\cT{\Cal T}

\def\tS{\tilde{S}}

\def\oN{\overline{N}}

\def\Cal{\mathcal}

\def\ov{\overline}

\def\Efi#1#2#3#4#5{\displaystyle
#1\!\!-\!\!#2
\!\!-\!\!#3
\!\!-\!\!#4
\hskip-24.2pt\lower4.5pt\hbox{${\scriptstyle|}
\hskip-3.35pt\lower6pt\hbox{$#5$}$}}

\def\Evia#1#2#3#4#5{\displaystyle
#1\!\!-\!\!#2
\!\!-\!\!#3
\hskip-24.2pt\lower4.5pt\hbox{${\scriptstyle|}
\hskip-3.35pt\lower6pt\hbox{$#4\!\!-\!\!\!-\!\!\!-\!\!$}$\hskip2.3pt$ 
{\scriptstyle|}
\hskip-3.35pt\lower6pt\hbox{$#5$}$}}

\def\Ezia#1#2#3#4{\displaystyle
#1\!\!-\!\!#2
\hskip-14.8pt\lower4.5pt\hbox{${\scriptstyle|}
\hskip-3.35pt\lower6pt\hbox{$#3\!\!-\!\!$}$\hskip2.3pt${\scriptstyle|}
\hskip-3.35pt\lower6pt\hbox{$#4$}$}}

\def\Efia#1#2#3#4#5#6{\displaystyle
#1\!\!-\!\!#2
\!\!-\!\!#3
\!\!-\!\!#4
\hskip-24.2pt\lower4.5pt\hbox{${\scriptstyle|}
\hskip-3.35pt\lower6pt\hbox{$#5$}$\hskip5.7pt${\scriptstyle|}
\hskip-3.35pt\lower6pt\hbox{$#6$}$}}

\def\Esi#1#2#3#4#5#6{\displaystyle
#1\!\!-\!\!#2
\!\!-\!\!#3
\!\!-\!\!#4\!\!-\!\!#5
\hskip-24.2pt\lower4.5pt\hbox{${\scriptstyle|}
\hskip-3.35pt\lower6pt\hbox{$#6$
\lower3pt\hbox{\ }}$}}

\def\Esia#1#2#3#4#5#6#7{\displaystyle
#1\!\!-\!\!#2
\!\!-\!\!#3
\!\!-\!\!#4\!\!-\!\!#5
\hskip-24.2pt\lower4.5pt\hbox{${\scriptstyle|}
\hskip-3.35pt\lower6pt\hbox{$#6$\hskip-3.8pt\lower4.5pt\hbox{$ 
{\scriptstyle|}
\hskip-3.35pt\lower6pt\hbox{$#7$}$}}
\lower3pt\hbox{\ }$}}

\def\Ese#1#2#3#4#5#6#7{\displaystyle
#1\!\!-\!\!#2
\!\!-\!\!#3
\!\!-\!\!#4\!\!-\!\!#5\!\!-\!\!#6
\hskip-33.6pt\lower4.5pt\hbox{${\scriptstyle|}
\hskip-3.35pt\lower6pt\hbox{$#7$
\lower3pt\hbox{\ }
}$}}

\def\Esea#1#2#3#4#5#6#7#8{\displaystyle
#1\!\!-\!\!#2
\!\!-\!\!#3
\!\!-\!\!#4\!\!-\!\!#5\!\!-\!\!#6\!\!-\!\!#7
\hskip-33.6pt\lower4.5pt\hbox{${\scriptstyle|}
\hskip-3.35pt\lower6pt\hbox{$#8$
\lower3pt\hbox{\ }
}$}}

\def\Eei#1#2#3#4#5#6#7#8{\displaystyle
#1\!\!-\!\!#2
\!\!-\!\!#3
\!\!-\!\!#4\!\!-\!\!#5\!\!-\!\!#6\!\!-\!\!#7
\hskip-43.2pt\lower4.5pt\hbox{${\scriptstyle|}
\hskip-3.35pt\lower6pt\hbox{$#8$
\lower3pt\hbox{\ }
}$}}

\def\Eeia#1#2#3#4#5#6#7#8#9{{\displaystyle
#1\!\!-\!\!#2
\!\!-\!\!#3
\!\!-\!\!#4\!\!-\!\!#5\!\!-\!\!#6\!\!-\!\!#7\!\!-\!\!#8
\hskip-52.2pt\lower4.5pt\hbox{${\scriptstyle|}
\hskip-3.35pt\lower6pt\hbox{$#9$
\lower3pt\hbox{\ }
}$}}}

\def\os{\overline{S}}

\def\tS{\tilde{S}}
\def\ts{\tS}
\def\ts7{\tilde{S}_7}

\def\bZ{\Bbb Z}
\def\bC{\Bbb C}

\def\bR{\Bbb R}

\def\tS{\tilde S}

\def\oS{\overline{S}}

\def\cT{\Cal T}

\def\tS{\tilde{S}}

\def\oN{\overline{N}}

\def\Cal{\mathcal}

\def\val{\operatorname{val}}
\def\oN{\overline{N}}
\def\os7p{\oS_7'}
\def\os7{\oS_7}
\def\on6{\oN_6}
\def\n6{\oN_6}

\newtheorem{Theorem}{Theorem}[section]

\newtheorem{Proposition}[Theorem]{Proposition}
\newtheorem{Lemma}[Theorem]{Lemma}

\newtheoremstyle{myreview}{}{}{}{}{\scshape}{.}{ }{}
\theoremstyle{myreview}

\theoremstyle{definition}

\newtheorem{Example}[Theorem]{Example}

\theoremstyle{remark}
\newtheorem{Remark}[Theorem]{Remark}

\newcounter{et}[Theorem]

\def\bZ{\mathbb Z}
\def\bC{\mathbb C}

\def\bR{\mathbb R}

\def\tS{\tilde S}

\def\oS{\overline{S}}
\def\cT{\Cal T}

\def\tS{\tilde{S}}

\def\oN{\overline{N}}

\def\Cal{\mathcal}

\def\ov{\overline}

\newcommand{\CC}{\mathbb C}
\newcommand{\RR}{\mathbb R}

\begin{document}

\title{The Newton Polytope of the Implicit Equation}

\author{Bernd Sturmfels}
\address{Department of Mathematics, University of California, Berkeley,
CA 94720-3840}
\email{bernd@math.berkeley.edu}

\author{Jenia Tevelev}
\address{Department of Mathematics,
  University of Massachusetts,
Amherst, MA 01003-9305}
\email{tevelev@math.utexas.edu}

\author{Josephine Yu}
\address{Department of Mathematics, University of California, Berkeley,
CA 94720-3840}
\email{jyu@math.berkeley.edu}
\thanks{Josephine Yu was supported by the NSF Graduate Research  
Fellowship.
Bernd Sturmfels was partially supported by NSF grant
DMS-0456960.}


%\date{\today }

\dedicatory{Dedicated to Askold Khovanskii on the occasion of his  
60th birthday}

\begin{abstract}
We apply tropical geometry to  study the image of a map defined by Laurent polynomials with generic coefficients. If this image is
a hypersurface then our approach gives
a construction of its Newton polytope.
\end{abstract}

\maketitle

\section{Introduction}

Implicitization is a fundamental operation
in computational algebraic geometry. Its objective is
to transform a given parametric representation of an algebraic variety
into its implicit representation as the zero set of polynomials.
Most algorithms for implicitization are  based on
multivariate resultants or Gr\"obner bases, but current implementations
of these algorithms tend to be  too slow to handle large instances.

Adopting Khovanskii's philosophy on the importance of Newton polytopes
\cite{Kho},
the foremost case of the implicitization problem is as follows.
Let $f_1,\ldots,f_n\in \CC[t_1^{\pm 1},\ldots,t_d^{\pm  1}]$
be Laurent polynomials with Newton polytopes
$P_1,\ldots,P_n \subset \mathbb{R}^d$ and with supports $A_i\subset  
P_i$.
We assume that each $f_i$ is {\em generic relative to its support},
which means that the coefficient vector of  $f_i$ lies in a
Zariski open  subset of $\CC^{A_i}$.
The aim of implicitization is to  compute the prime ideal $I$ of
all polynomials $g \in \CC[x_1,\ldots,x_n]$ which satisfies
$g(f_1,\ldots,f_n) \equiv 0$ in
$\CC[t_1^{\pm 1},\ldots,t_d^{\pm  1}]$.
We here seek to read off
as much information as possible about
the variety $V(I)$ from $P_1,\ldots,P_n$.

Of particular interest is the case $n=d+1$,
when $I = \langle\, g\, \rangle$ is a principal ideal.
Here the problem is to predict the
Newton polytope $Q$ of the hypersurface
$\,V(I) = \{ g = 0\}\,$ from $P_1,\ldots,P_n$.
This problem was posed
in \cite{SY} and has remained open for over a decade.
It also reappeared in recent work of  Emiris and Kotsireas \cite{EK}.

We here present a general solution to the problem of \cite{SY},
namely, a construction of the
Newton polytope $Q$ in terms of the input polytopes $P_i$.
When the variety $V(I)$ is not a hypersurface
but has codimension greater than one, the role of $Q$ is played by the
Chow polytope \cite{KSZ}, and our construction generalizes to that
case.

The present exposition is not self-contained;
it refers for two proofs to forthcoming papers.
Our~calculations are based on geometric characterizations of $\mathcal 
{T}(I)$
derived by Hacking, Keel and Tevelev in \cite{HKT, Te}.
A theory of {\em Tropical Implicitization}, which further
develops these results, appears in \cite{STY}.
In that theory, the coefficients of the $f_i(t)$ can  be arbitrary  (non-generic)
complex numbers, or even scalars
in an arbitrary field with a non-archimedean valuation.
One objective of this paper
is to announce and explain tropical implicitization
to a wider audience.

Let us begin by examining the  simplest
case: a parametrized curve in the plane.

\begin{figure}[htb]
\includegraphics[scale=0.9]{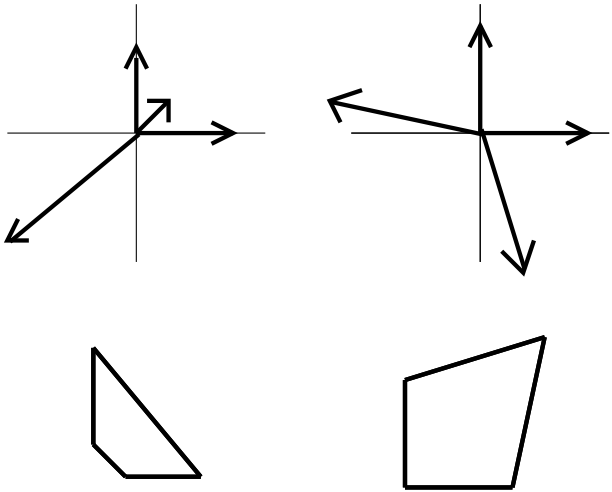}
\caption{Tropical plane curves and their Newton polygons}
\label{fig:planeCurves}
\end{figure}

\begin{Example}
\label{CurveInThePlane0}
Let $d=1$ and $n=2$. We wish to compute the
equation $g(x_1,x_2)$ of the plane algebraic curve
parametrized by two Laurent polynomials
$ \,x_1 = f_1(t)$ and $x_2 = f_2(t)$.
The Newton polytopes of $f_1$ and $f_2$
are segments on the  line $\RR^1$,
\begin{equation}
\label{segments}
P_1 \,=\, [\,a \,,\, b\,]  \quad
  \hbox{and} \quad P_2 \,=\, [\,c\, ,\,d\,] .
\end{equation}
Here $a \leq b$ and $c \leq d$ are integers.
The Newton polygon $Q$ of
the implicit equation  $g(x_1,x_2)$ is
  a (possibly degenerate)  quadrangle in  $\RR^2$.
There are four cases:
\hfill \break
$\bullet$ if $ a \geq 0 $ and $c \geq 0 $ then $\, Q = {\rm conv}
\bigl\{ (0, b), (0, a), (c, 0), (d, 0) \bigr\}$,
\hfill \break
$\bullet$ if $b \leq 0  $ and $d \leq 0 $ then $\,  Q = {\rm conv}  
\bigl\{
  (0, -a), (0, -b), (-d, 0), (-c, 0)  \bigr\}$,
\hfill \break
$\bullet$ if $ a \leq 0 $, $d \geq 0 $
and $bc \geq ad$ then $\,  Q = {\rm conv}
\bigl\{  (0, b-a), (0, 0), (d-c, 0), (d, -a) \bigr\}$,
\hfill \break
$\bullet$ if $ b \geq 0 $, $c \leq 0 $
and $bc \leq ad$ then $\,  Q = {\rm conv} \bigl\{
(0, b-a), (0, 0), (d-c, 0), (-c, b)  \bigr\}$.
\hfill \break
The possible quadrangles $Q$ and their normal fans
are illustrated in Figure \ref{fig:planeCurves}.
We note that  the case $a=48, b=63, c=d=32$
appears in \cite[Example 3.3]{EK}.
\qed
\end{Example}

Tropical implicitization can be visualized using Figure \ref 
{fig:planeCurves}.
Given the input  (\ref{segments}), what we compute is the
tropical curve $\mathcal{T}(g)$. This is
the picture in the first row. By definition,
$\mathcal{T}(g)$ is the set of rays in the
inner normal fan of the (unknown) polygon in the second row.
The output is
the collection  of four rays which are spanned
by the vectors $(b-a,0)$, $(0,d-c)$,
$(a,c)$,  and  $(-b,-d)$.  These vectors sum to zero,
and they determine the polygon $Q$ according
to the four cases spelled out above.

The motivation behind both \cite{EK} and \cite{SY} 
is that {\it a priori} knowledge of the Newton polytope $Q$
would greatly facilitate the subsequent computation
of recovering the coefficients of $g(x)$
from the coeffcients of $f_1(t),\ldots,f_n(t)$.
This is a problem of numerical linear algebra,
and it will be discussed in Subsection \ref{subsec:CI}.

Our presentation is organized into four sections as follows.
In Section 2 we present our first main result,
Theorem \ref{SolutionToYuProblem}, which is a formula
for the tropical variety $\mathcal{T}(I)$
in terms of the Newton polytopes $P_1,\ldots,P_n$.
In Section 3 we demonstrate applications of this
formula for a wide range of examples.
In Section 4 we give a formula, in terms of mixed volumes and lattice  
indices, for the intrinsic multiplicities
of the maximal cones in the tropical variety $\mathcal{T}(I)$.
Knowledge of these multiplicities is essential when trying
to reconstruct information about the ideal $I$ from
its tropical variety $\mathcal{T}(I)$.
In Section 5 we present our algorithm for computing the
Newton polytope $Q$ of the implicit equation
(or the Chow polytope of $I$)  from $P_1,\ldots,P_n$.

\section{The tropical variety of the implicit equations}

Our implicitization problem is
specified by a collection of $n$ Laurent polynomials
\begin{equation}
\label{InputData}
  f_i (t) \quad = \quad \sum_{a \in A_i} c_{i,a} \cdot t_1^{a_1}  
\cdots t_d^{a_d}
\qquad \qquad  (i =1,2,\ldots,n) .
\end{equation}
Here each $A_i$ is a finite subset of $\bZ^d$, and the $c_{i,a}$
are generic complex numbers.
Our ultimate aim is to compute the ideal
$\,I \subset \CC[x_1,\ldots,x_n]\,$  of
algebraic relations among $f_1(t), \ldots,f_n(t)$,
or at least, some information about its variety $V(I)$.

The tropical approach to this problem
is based on the following idea. Rather than computing
$V(I)$ by algebraic means,
we shall compute the
tropical variety $\mathcal{T}(I)\subset\bR^n$ by combinatorial means.
We recall (e.g.~from \cite{BJSST, PS, SS}) that
the {\em tropical variety} $\mathcal{T}(I)$ is the set
of all vectors $w \in \bR^n$ such that
the {\em initial ideal} $\,{\rm in}_w(I) = \langle
{\rm in}_w(g) : g \in I \rangle \,$ contains no monomials.
Here ${\rm in}_w(g)$ is the {\em initial form} of $g$ which is
the sum of all terms $\,c_\alpha x^\alpha\,$
in $g$ which have $w$-minimal weight $\, \alpha \cdot w $.

Let  $\Psi : \RR^d \rightarrow \RR^n$ be the tropicalization of the  
map $f = (f_1, \dots, f_n)$, i.e.,
let $ \Psi_i(w) = \text{min}\{w \cdot v : v \in P_i\} $ be the  
support function of the Newton polytope $P_i = {\rm conv}(A_i)$.
The image of $\Psi$ is contained in the
tropical variety $\mathcal{T}(I)$, by~\cite[Theorem 2]{PS},
but this containment is usually strict.
In other words, the image of the tropicalization of $f$ is usually a proper  
subset of the tropicalization of the image of $f$.
The point of the following result is to
characterize the difference  $\mathcal{T}(I)\backslash {\rm image} 
(\Psi)$.

Let $e_1,\ldots,e_n$ be the standard basis of $\bR^n$. For
  $J \subseteq \{1, \dots, n\}$, we abbreviate the orthant
$\, \mathbb{R}_{\geq 0} \{e_j: j \in J\}\,$ by
  $\,\mathbb{R}_{\geq 0}^J \,$
  and the Minkowski sum
  $\, \sum_{j \in J} P_j\,$ by $\,P_J  $.

\begin{Theorem}
\label{SolutionToYuProblem}
Let $f_1, \ldots,f_n$ be Laurent polynomials
in $d$ variables which are generic
relative to their support, and let
$\,I \subset \CC[x_1,\ldots,x_n]\,$ be the
prime ideal of all algebraic relations among
$f_1,\ldots,f_n$.
Then the following subsets of $\mathbb{R}^n$ coincide:
\begin{enumerate}
\item the tropical variety $\,\mathcal{T}(I)$,
\item the union of all sets
  $\,\Psi (\mathcal{T}(\langle f_j : j \in J \rangle)) \,+ \,
\mathbb{R}_{\geq 0}^J $,
where $\,J \subseteq \{1,2,\ldots,n\}$,
\item the union of the cones $\,\Psi(w) \,+  \,  \mathbb{R}_{\geq 0} 
^J $,
where $w \in \mathbb{R}^d$ and $J$ runs over all subsets of $\{1,2, 
\ldots,n\}$ such that
$\,\langle\, \text{in}_w(f_j) : j \in J \,\rangle\,$ contains no  
monomial,
\item
  the union of the cones $\,\Psi(w) \,+  \,  \mathbb{R}_{\geq 0}^J\,$
  such that, for all subsets $K \subseteq J$, the linear form $w$
attains its minimum over  $P_K$ at a face of
  dimension $ \geq |K|$.
  \end{enumerate}
\end{Theorem}

Characterization (4) gives rise to a
combinatorial algorithm for computing
the tropical variety $\mathcal{T}(I)$
  directly from the given Newton
polytopes $P_1,\ldots,P_n$.
We will see examples in Section 3, and in Section 5
we discuss our implementation of this algorithm.
In this section we explain Theorem  \ref{SolutionToYuProblem} and
give an outline of its proof.

It is instructive to note that the contribution of
the empty set $J = \emptyset$ in Theorem \ref{SolutionToYuProblem}  
(2) is precisely the image of the tropicalization $\Psi$ of the given  
map $f$:
\begin{equation}
\label{ImageOfPsi}
  \Psi (\mathcal{T}(\langle  \emptyset \rangle)) \,+ \,
\mathbb{R}_{\geq 0}^\emptyset
\,\,=\,\,  \Psi (\mathcal{T}(\{0\})) \,\, = \,\,
\Psi(\RR^d) \,\, = \,\, {\rm image}(\Psi).
\end{equation}
Thus it is the contributions made by
non-empty subsets $J$ which make up the
difference between the tropicalization of
the image and the image of the tropicalization.

\begin{Example}
\label{CurveInThePlane2}
Consider the case of
a plane curve as in Example \ref{CurveInThePlane0}.  Assume for  
simplicity that $f_1$ and $f_2$ are not monomials, i.e.,
$a < b$ and $c < d $.
For $J = \{1,2\}$ we just get the empty set
since $\,\langle f_1,f_2 \rangle \,$ is the unit
ideal in $\CC[t,t^{-1}]$.
For $J = \{1\}$ we get the ray spanned by $(1,0)$
since  $\,\cT(\langle f_1 \rangle) = \{0\}$, and for
$J = \{2\}$ we get the ray spanned by $(0,1)$.
Finally, for $J = \emptyset$ we get the image of  the tropicalization
$\Psi$ which consists of the two rays spanned by
$(a,c)$ and by $(-b,-d)$. See Figure \ref{fig:planeCurves}.
\qed
\end{Example}

We continue by explaining the meaning of the combinatorial conditions  
(2), (3) and (4),
and, in the course of doing so, we shall prove that they are equivalent.
First of all, we relabel so that $J = \{1,\ldots,r\}$
and consider the system of equations
\begin{equation}
\label{initialsystem}
  {\rm in}_w(f_1)(t)\, =\, \cdots \,=\, {\rm in}_w(f_r)(t)\, =\, 0 .
  \end{equation}
   The Newton polytopes of these Laurent polynomials
  are $\, {\rm face}_w(P_1) , \ldots, {\rm face}_w(P_r)$,
  where ${\rm face}_w(P_i)$ is the face of $P_i$
  at which the linear functional $u \mapsto u \cdot w$ attains
  its minimum. Let $s$ denote the dimension of the polytope
  \begin{equation}
\label{facewPJ}
  {\rm face}_w(P_J) \, = \,  {\rm face}_w(P_1) + \cdots +  {\rm face}_w(P_r)
  \end{equation}
  We wish to determine whether or not the system
  (\ref{initialsystem}) has a solution $t$   in the algebraic torus  $ 
(\bC^*)^d$.
  The equations (\ref{initialsystem}) can be rewritten as a system of
  Laurent polynomials in $s$ unknowns. Since the coefficients of the
  $f_i$ are assumed to be generic, a necessary condition for (\ref 
{initialsystem})
  to be solvable is that $r \leq s$. Let us first consider the case  
$r = s$.
   By Bernstein's Theorem \cite{Be, Kho}, the number of solutions
  to the system (\ref{initialsystem}) equals
  \begin{equation}
  \label{MV}
   {\rm MixedVolume}\bigl( {\rm face}_w(P_1) , \ldots, {\rm face}_w(P_r) \bigr) ,
   \end{equation}
  where the mixed volume is normalized with respect to the
  lattice parallel to the affine span of ${\rm face}_w(P_J)$.
In degenerate cases this mixed volume may be zero.

\begin{Lemma}
\label{MVpositive}
Suppose $r = s$. Then the following conditions are equivalent:
\begin{enumerate}
\item The mixed volume (\ref{MV}) is positive.
\item The system (\ref{initialsystem}) has at least one solution $t  
\in (\bC^*)^d$.
\item For all subsets $K \subseteq J$, the polytope
  ${\rm face}_w(P_K)$ has
  dimension $ \geq |K|$.
\end{enumerate}
\end{Lemma}

\begin{proof}
The equivalence of (1) and (2) follows from Bernstein's Theorem \cite 
{Be}.
The equivalence of (1) and (3) appears in
  \cite[Theorem 4.13, page 127]{Ew}.
\end{proof}

Next suppose that $r < s$. Then we can artificially
add $s-r$ equations with generic coefficients to obtain a system
of $s$ equations in $s$ unknowns. Application of Lemma \ref{MVpositive} to  
this new system shows that the dimension criterion continues to hold.
Furthermore, the following proposition shows that the solvability of  
(\ref{initialsystem}) is equivalent to the
conditions appearing  in Theorem \ref{SolutionToYuProblem}.

\begin{Proposition}
\label{fourequiv}
For any $w \in \bR^d$ and $J = \{1,\ldots,r\}$
the following are equivalent:
\begin{enumerate}
\item The initial system (\ref{initialsystem}) has at least one  
solution $t \in (\bC^*)^d$.
\item For all subsets $K \subseteq J$, the polytope
  ${\rm face}_w(P_K)$ has
  dimension at least $ |K|$.
  \item The  ideal
  $\,\langle\, \text{in}_w(f_j) : j \in J \,\rangle\,$ contains no  
monomial.
\item The vector $w$ lies in $\, \mathcal{T}(\langle f_j : j \in J  
\rangle) $.
\end{enumerate}
\end{Proposition}

\begin{proof}
The equivalence of (1) and (2) follows from Lemma \ref{MVpositive}.
Since monomials are units in the Laurent polynomial ring,
Hilbert's Nullstellensatz shows that (1) and (3) are equivalent.
Now (4) holds if and only if there is no monomial in
$\,{\rm in}_w \bigl(\langle f_1, \ldots,f_r \rangle \bigr) \,$.
This condition implies (3) but, a priori, it may be stronger.
To see that they are  equivalent, we will use the hypothesis
that the coefficients of the $f_i$ are generic and show that (1)  
implies (4).
Suppose that (1) holds, and let $t^0 \in (\bC^*)^d$ be a solution of  
the system of equations  (\ref{initialsystem}).
Using the polyhedral homotopy of \cite{HS}, we can
construct a solution $t^1$ of the equations
$\,f_1(t) = \cdots = f_r(t) = 0$ which degenerates
to $t^0$ under the one-parameter torus given by $w$.
The deformation $t^0 \rightarrow t^1$
ensures that $t^0$ is a point in the variety of the initial ideal
$\,{\rm in}_w \bigl(\langle f_1,\ldots,f_r \rangle \bigr) $,
which therefore contains no monomial.
\end{proof}

\begin{Remark}
Our deformation argument for \ (1) $\Longrightarrow$ (4)\  can be
understood geometrically using the
smooth projective toric variety $Y$ below.
Namely, both $t^0$ and $t^1$ are points on $Y$ that lie in the  
intersection of divisors $\ov E_1 \cap \cdots \cap \ov E_r $
but $t^1$ lies in the dense torus while $t^0$ lies in the toric  
boundary in direction $w$.
\end{Remark}

Proposition
\ref{fourequiv} shows that the last three conditions in Theorem
\ref{SolutionToYuProblem} are equivalent. To make the connection
to the first condition, we need a few ingredients from algebraic  
geometry. Let
$E_i=\{t\in  (\CC^d)^*: f_i(t)=0\}\,\,\text{and}\,\,Y^0= (\CC^*)^d  
\setminus\bigcup_{i=1}^nE_i.$
The Laurent polynomials (\ref{InputData}) specify  a morphism
of affine algebraic varieties
\begin{equation}\label{mapF}
f:\,Y^0\to(\bC^*)^n,\quad z\mapsto(f_1(z),\ldots,f_n(z))
\end{equation}
Our goal is to compute the tropicalization of the image of $f$.
Suppose that  $Y \supset Y^0$ is any smooth compactification of $Y^0$
whose boundary $D = Y \backslash Y^0$ is a
divisor with simple normal crossings.  Let $D_1,\ldots,D_m$ denote
the irreducible components
of the boundary divisor $D$. Let $ \Delta_{Y,D}$ denote the
  simplicial complex on $\{1,\ldots,m\}$ whose simplices are the
  indices of subsets of divisors $D_i$ whose intersection on $Y$
  is non-empty. Since $D$ is a normal crossing
  divisor, the simplicial complex $\Delta_{Y,D}$ is pure of dimension  
$d-1$.
For any of the irreducible divisors $D_i$, let $[D_i]\in\bZ^n$ be the  
vector $(\val_{D_i}f_1,\ldots,\val_{D_i}f_n)$, where $\val_{D_i}$ is  
the order of zeros-poles along $D_i$.

\begin{Theorem} \label{tropforncc}
The tropical variety $\cT(I)$ is equal to the union of cones
\begin{equation}
\label{unionofcones}
  \cT(I) \,\,\,= \,\,\, \bigcup_{\sigma \in \Delta_{Y,D}}
\RR_{\geq 0} \bigl\{ [D_i] \,: \, i \in \sigma \bigr\} .
\end{equation}
\end{Theorem}

The proof of this theorem will appear in \cite{STY}
and follows the work in  \cite{HKT}.
There, $Y^0$ may be
any variety that is {\em very affine}, which means that it is closed
in a torus, and $f$ is allowed to be any morphism of
very affine varieties. The required compactification $Y$ of
the very affine variety $Y^0$
can be constructed, in principle, by appealing to Hironaka's
Theorem on resolution of singularities. In practice, however,
we are interested in situations where $Y$ is more
readily available.

This is precisely the case studied in this paper,
namely, when the coefficients $c_{i,a}$ of the Laurent
polynomials $f_i(t)$ are assumed to be generic.
In this case the desired compactification $Y$ of $Y^0$
is furnished to us by the methods of toric geometry.

\begin{proof}[Proof of Theorem~\ref{SolutionToYuProblem}]
The equivalence of the sets in (2), (3) and (4)
was shown in Proposition \ref{fourequiv}. What remains to be shown is  
the
equivalence of (1) and (2).

Let $Y$ be a $d$-dimensional smooth projective toric variety whose
polytope $\tilde{P}$ has the given Newton polytopes $P_1,\ldots,P_n$
as  Minkowski summands. The smooth toric variety $Y$ is
a compactification of $Y^0$. It comes with a canonical morphism
$Y \rightarrow Y_i$ onto the (generally not smooth)  toric variety $Y_i$
associated with $P_i$.

We claim that its
boundary $Y \backslash Y^0$ has simple normal crossings.
The irreducible components of $Y\backslash Y^0$ are of two types.
Firstly, we have toric divisors $D_1,\ldots,D_l$ indexed by the  
facets of $\tilde{P}$.
The toric boundary $\,D_1 \cup \cdots \cup D_l\,$ of  $Y$ has simple  
normal crossings because
$\tilde{P}$ is a simple smooth polytope.
Secondly, we have divisors $\ov E_1,\ldots,\ov E_n$ which
are the closures in $Y$ of the divisors $E_1,\ldots,E_n$ in $(\bC^*)^d$.
The divisor $\ov E_i$ is the pullback under the morphism $Y  
\rightarrow Y_i$
of a general hyperplane section of the projective embedding of $Y_i$
defined by $P_i$. Bertini's Theorem implies that
the $\ov E_i$ are smooth and irreducible, and that
the union of all $D_i$'s and all $\ov E_j$'s has normal crossings.
(Here we are tacitly assuming that each polytope $P_i$
has dimension $\geq 2$. If ${\rm dim}(P_i) = 1$ then
$\ov E_i $ is the disjoint union of smooth and irreducible
divisors, and the following argument needs to be slightly modified).

In summary, we have that Theorem \ref{tropforncc} can be applied to
$$ Y \backslash Y^0 \quad = \quad
  D_1 \cup D_2 \cup \cdots \cup D_l\,\cup\,
  \ov E_1 \cup \cdots \cup \ov E_n .  $$
The simplicial complex $\Delta_{Y,D}$ has dimension $d-1$
and it has $m = l+n$ vertices, one for each of the divisors
$D_i$ and $\ov E_j$. Its maximal simplices correspond to pairs
$(C,J)$ where $C = \{i_1,\ldots,i_{d-r} \} \subseteq \{1,\ldots,l\}$ and
$J  =\{j_1,\ldots,j_r\}\subseteq \{1,\ldots,n\}$ and
\begin{equation}
\label{intersectonY}
D_{i_1} \cap \cdots \cap D_{i_{d-r}} \,\cap\, \ov E_{j_1} \cap \cdots  
\cap \ov E_{j_r} \quad
\not= \,\,\,\, \emptyset.
\end{equation}
For any $\,J \subseteq \{1,\ldots,n\}\,$ let $\Delta_{Y,D}[J]$
denote the subset of $\Delta_{Y,D}$ consisting of all simplices
with fixed $J$. Note that $\Delta_{Y,D}[\emptyset]$
is the boundary complex
of the simplicial polytope dual to $\tilde{P}$. Moreover,
$\Delta_{Y,D}[J]  = \{ \emptyset \}$ if $|J|= d$,
and $\Delta_{Y,D}[J]  = \emptyset$ if $|J| > d$.
For each $j \in J$ we have
$\,[\ov E_j] = (\val_{\ov E_j}f_1,\ldots,\val_{\ov E_j}f_n) \, = \, e_j \,$
is the $j^{\rm th}$ basis vector in $\bR^n$. Hence
the formula (\ref{unionofcones}) in Theorem \ref{tropforncc} can be
rewritten as follows:
$$  \mathcal{T}(I) \quad = \quad
\bigcup_{J \subset \{1,\ldots,n\}}
\bigl(\,\, \bR_{\geq 0}^J \,\,+ \!\!
\bigcup_{C \in \Delta_{Y,D}[J]}
\RR_{\geq 0} \bigl\{ [D_i] \,: \, i \in C \bigr\} \, \bigr). $$
  Hence to prove
  the remaining equivalence (1) = (2), it suffices to show that
\begin{equation}
\label{remainstoshow}
\Psi (\mathcal{T}(\langle f_j : j \in J \rangle)) \quad = \quad
\bigcup_{C \in \Delta_{Y,D}[J]}
\RR_{\geq 0} \bigl\{ [D_i] \,: \, i \in C \bigr\}  .
\end{equation}
Let $w_i \in \bR^d$ denote the primitive inner normal vector of
the facet of the polytope $\tilde{P}$ corresponding to the divisor $D_i$.
For any $j \in J$, the integer ${\rm val}_{D_i}(f_j)$
is the  value at $w_i$ of the support function of the polytope $P_j$.
  Hence $[D_i] = \Psi(w_i)$ in $\bR^n$, so right hand side of (\ref{remainstoshow}) is the image under $\Psi$ of the subfan of the normal fan of $\tilde{P}$ indexed by
$\Delta_{Y,D}[J]$.  Thus it suffices to show that this subfan coincides with $\cT(\langle f_j : j \in J \rangle)$.  A vector $w$ lies in this subfan if and only if (\ref{intersectonY}) holds, which is the same as saying that the system ${\rm in}_w(f_j) = 0$ for $j \in J$ has a solution in $(\CC^*)^d$.  By Proposition \ref{fourequiv}, this is equivalent to $w \in \mathcal{T}(\langle f_j : j \in J \rangle)$.
This completes our proof of Theorem \ref{SolutionToYuProblem}.
\end{proof}

\section{Examples}

We now present five families of examples which
illustrate the Theorem \ref{SolutionToYuProblem}.

\subsection{Monomials}
Suppose that each $f_i$ is a monomial, so each $P_i$ is just a point,  
say
$P_i = \{a_i\} \subset \bZ^d$.
Then $I$ is the toric ideal generated by all binomials
$x^u - x^v$ where $u-v$ is in the kernel
of the matrix $A = [a_1,a_2,\ldots,a_n]$.
The tropical variety $\mathcal{T}(I)$ equals
the row space of $A$.
The tropical morphism $\Psi : \bR^d \rightarrow \bR^n$
is the linear map $w \mapsto w A$ given by left multiplication with $A$. For any non- 
empty subset $J$
of $\{1,2,\ldots,n\}$,  the ideal $\langle f_j : j \in J \rangle$ is  
the unit ideal,
which contributes nothing to the union
in (2) of Theorem \ref{SolutionToYuProblem}.
The ideal corresponding to  $J = \emptyset$ is
$\langle \emptyset \rangle = \{0\}$, and
$\,\mathcal{T}(I) = {\rm rowspace}(A)\,$ is indeed the image of
$\,\mathcal{T}(\langle \emptyset \rangle) = \bR^d\,$
under the map $\Psi$.

\subsection{The Unmixed Case}
Suppose that all Newton polytopes are equal,
say $\,P_1=\ldots=P_n=P$, and $P$ has the maximal
dimension $d$.  We distinguish two cases.

First suppose that $P$ contains the origin.
Then the image of $\Psi$ is the halfline
$\,C^-=\bR_{\le0}(e_1+\ldots+e_n)$. In fact,
the set $\,\Psi (\mathcal{T}(\langle f_j : j \in J \rangle))\,$
in Theorem \ref{SolutionToYuProblem} (2)
equals $\,C^-\,$ if $\,|J| < d$, it equals
$\,\{0\}\,$ if $\,|J|=d$, and it is empty if $\,|J| > d$.
Therefore, the tropical variety
$\cT(I)$ is the union of the cones
$\, \bR_{\geq 0}^J$ for $J\subset\{1,\ldots,n\}$ with $|J|=d$
and the cones $\,C^- + \bR_{\geq 0}^J\,$
for $J\subset\{1,\ldots,n\}$ with $|J|=d-1$.

Next suppose that $P$ does not contain the origin.
Then the image of $\Psi$ is the line
$\,C =\bR(e_1+\ldots+e_n)$, and we conclude that
the tropical variety
$\,\cT(I) \,$ is the union of the cones
$ \,\bR_{\geq 0}^J\,$ for $\,|J|=d$,
and the cones $\,C + \bR_{\geq 0}^J \,$ for  $\,|J|=d-1$.

In the case $n=d+1$, the tropical variety
$\cT(I)$ just constructed has codimension
one in $\bR^n$. It is the union of all codimension one
cones in the normal fan of the Newton
polytope $Q$ of the implicit equation. The polytope $Q$
was constructed with a different technique in
\cite[Theorem 9]{SY}, and the previous two paragraphs are
consistent with that result. However, knowledge
of $\cT(I)$ is not sufficient to determine $Q$.
For that we need the formula for multiplicities
in Section 4 below.

\begin{Example}
Let $P$ be any convex lattice polygon in $\bR^2$
and consider three generic bivariate Laurent
polynomials with Newton polygon $P$. Then
$\cT(I)$ is a tropical surface in $\bR^3$.
We abbreviate $e_4 = -e_1-e_2-e_3$.
If the polygon $P$ contains the origin in $\bR^2$,
then $\cT(I)$ is the fan consisting
of the six two-dimensional cones
$$ \bR_{\geq 0} \{ e_1,e_2 \}, \,
\bR_{\geq 0}\{ e_1,e_3 \},\,
\bR_{\geq 0}\{ e_2,e_3 \},\,
\bR_{\geq 0}\{ e_1,e_4 \} , \,
\bR_{\geq 0}\{ e_2,e_4\} , \,
\bR_{\geq 0}\{ e_3,e_4 \}. $$
This implies  that $Q$ is a tetrahedron
$\,{\rm conv}\{0, \alpha e_1, \alpha e_2, \alpha e_3\}$,
for some $\alpha > 0$.

If $P$ does not contain the origin,
then the tropical surface
  $\cT(I)$ has a natural fan  structure
which is given by the
nine two-dimensional cones
$$ \bR_{\geq 0} \{ e_i,e_j \}, \,
\bR_{\geq 0}\{ e_i,e_4 \} , \,
\bR_{\geq 0}\{ e_i,-e_4 \} \qquad \qquad (1 \leq i < j \leq 3). $$
This implies that $Q$ is a triangular prism
$\,{\rm conv}\{
\alpha e_1, \alpha e_2, \alpha e_3,
\beta e_1, \beta e_2, \beta e_3\}$,
for some $\alpha > \beta > 0$.
The parameters $\alpha, \beta$ are
expressed in terms of volumes as in
\cite[Theorem 9]{SY}. Alternatively, we can
determine $\alpha$ and $\beta$ using
Theorem  \ref{multFormula} below. \qed
\end{Example}

\subsection{Linear maps}
Suppose each support $A_i$ is a subset of
$\{e_1,\ldots,e_d\}$, the unit vectors in $\bR^d$.
Thus $f_i$ is a generic linear form in the set of unknowns
$\,\{\,t_j \,: \, e_j \in A_i\}$. The ideal $I$ is generated
by an $(n-d)$-dimensional space of linear forms in
$\,\bC[x_1,\ldots,x_n]$.
We introduce the subset
$\,B_j = \{i : e_j \in A_i\}\,$ for $j = 1,2,\ldots,d$.
The rank $d$ matroid corresponding to the linear space
$V(I) \subset \bC^n$ is the {\em transversal matroid} \cite{Br}
of the set family $\{B_1,\ldots,B_d \}$.
The tropical variety $\mathcal{T}(I)$ is the
{\em Bergman complex} of that transversal matroid,
as shown in \cite{AK}.

Thus, in the linear case, Theorem \ref{SolutionToYuProblem}
offers an interesting new representation of the Bergman fans
of transversal matroids. It can be described as follows.
The sum of simplices $\,P = \sum_{i=1}^n {\rm conv}(A_i)\,$
is a {\em generalized permutohedron} \cite{Po}, which means
that the normal fan of $P$
is a coarsening of the $S_n$-arrangement. Let $C$ be any cone
in that fan. Each initial form
$\,{\rm in}_C(f_j)\,$ is supported on a subset
${\rm face}_C(A_j)$ of $A_j$.  Then with $C$ we associate the simplicial
complex $\Delta_C$ on $\{1,2,\ldots,n\}$ whose simplices are the subsets
$J$ such that these initial forms have a solution in $(\bC^*)^d$.
This is a combinatorial condition on the set family
$\,\{ \,{\rm face}_C(A_j) : j \in J \, \} $, which is essentially
Postnikov's {\em Dragon Marriage Condition} \cite[\S 5]{Po}.
To construct the Bergman fan, we map the normal fan
of $P$ from $\bR^d$ into $\bR^n$ using the map $\Psi$,
and to each image cone $\Psi(C)$ we attach the
family of orthants $\,\bR_{\geq 0}^J\,$ indexed by $\Delta_C$.

\begin{Example}
\label{trans}

\begin{figure}[htb]
\includegraphics[scale=1.0]{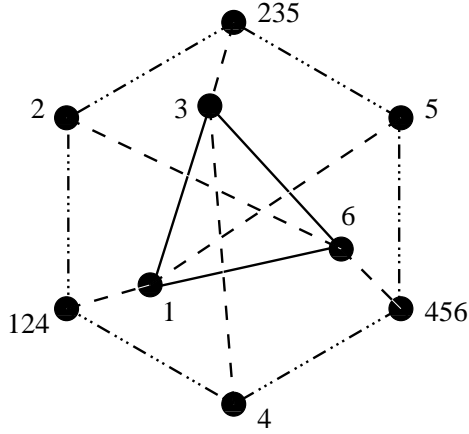}
\caption{The Bergman complex of a transversal matroid}
\label{fig:berg36}
\end{figure}

Let $d=3$, $n=6$ and consider six linear forms
$$
\begin{matrix}
f_1 &=&  c_{11} \cdot t_1\,+\, c_{12} \cdot t_2 \\
f_2 &=& c_{21} \cdot t_1 \\
f_3 & =& c_{31} \cdot t_1 \,+\, c_{33} \cdot t_3 \\
f_4 &=& c_{42} \cdot t_2 \\
f_5 &=& c_{53} \cdot t_3 \\
f_6 &= & c_{62} \cdot t_2 \,+\, c_{63} \cdot t_3,
\end{matrix}
$$
where the $c_{ij}$ are general complex numbers.
The tropical linear space $\cT(I)$ is three-dimensional in
$\bR^6$, but each maximal cone contains the line
spanned by $(1,1,1,1,1,1)$, so we can represent
$\cT(I)$ by a graph. This graph is the Bergman complex of a
rank $3$ matroid on $\{1,2,3,4,5,6\}$, namely, the transversal  
matroid of
$$ (B_1,B_2,B_3) \quad = \quad \bigl( \{1,2,3\},\,\{1,4,6\},\,\{3,5,6 
\} \bigr). $$

The drawing of $\cT(I)$ in Figure \ref{fig:berg36} is the same as   
\cite[Figure 2]{DFS},
but we now derive it from Theorem \ref{SolutionToYuProblem} (2).
The polytope $\,P = P_1 + \cdots +P_6\,$ is  a planar hexagon,
and $\Psi$ maps the six
two-dimensional cones in its normal fan to the six edges in the
outer hexagon in Figure \ref{fig:berg36}.  Each of the edges
$\{124,1\}, \{2,6\}, \{235,3\},\{5,1\}, \{456,6\}$, and $\{4,3\}$ are  
between the image of a ray in the normal fan of $P$ and a ray  
generated by $e_1,e_3$, or $e_6$.  For instance, the edge $\{124,1\}  
$ corresponds to the ray generated by $\Psi\left((1,1,0)\right) =  
(1,1,0,1,0,0)$ plus the ray $\bR_{\geq 0} e_1$ in $\bR^6$.
None of the simplicial complexes $\Delta_C$ contains $2, 4$, or $5$
because $f_2, f_4$, and $f_5$ are monomials.
  However, $e_2, e_4$, and $e_5$ are images of some rays in the  
normal fan of $P$ under the tropical map $\Psi$.
The triangle in the middle of Figure \ref{fig:berg36}
represents the simplicial complex associated with the
$0$-dimensional face $\{0\}$ of the normal fan:
$$ \Delta_{\{0\}} \,\,\, = \,\,\,
\bigl\{ \{1,3\}, \{1,6\}, \{3,6\} \bigr\}. $$
This shows how tropical implicitization works for the
linear forms $\,f_1,\ldots,f_6$.  \qed
\end{Example}

\subsection{Binomials}
Consider the case where the map $f$ is given by $n$ binomials
$$ f_i \,\,=\,\, c_{i1} \cdot t^{a_i} \,+ \, c_{i2} \cdot  t^{b_i}
\qquad (i=1,2,\ldots,n). $$
Each Newton polytope $\,P_i = {\rm conv}(a_i,b_i)\,$ is a line segment,
and their Minkowski sum $\,P = P_1 + \cdots + P_n\,$ is a {\em zonotope}
(i.e.~a projection of the $n$-cube).
The normal fan of $P$ is the {\em hyperplane arrangement}
$\,\mathcal{H} = \bigl\{
\{ u \cdot a_j = u \cdot b_j \}\bigr\}_{j=1,\ldots,n}$.
The map $\Psi : \bR^d \rightarrow \bR^n$
is the tropical morphism associated
with the arrangement $\mathcal{H}$, that is,
$$
\Psi(u) \,\,=\,\,
\bigl(
{\rm min}(u \cdot a_1, u\cdot b_1),
{\rm min}(u \cdot a_2, u\cdot b_2), \ldots,
{\rm min}(u \cdot a_n, u\cdot b_n) \bigr).
$$
This map was recently studied by Ardila \cite{Ar}
for certain graphic arrangements $\mathcal{H}$.

Theorem \ref{SolutionToYuProblem} shows that
the tropical variety $\cT(I)$ is the union of the cones
$\,\Psi(C) + \RR_{\geq 0}^J$, where $C$ is any cone
of the hyperplane arrangement $\mathcal{H}$, and
the hyperplanes indexed by  $\,J \subset \{1,\ldots,n\}$
contain the cone $C$ and
are linearly independent.

\begin{Example}
Let $d=2$ and consider the following $n=3$ binomials in $t_1$ and $t_2$:
$$ f_1 \,\,=\,\, (t_1 - \alpha) \cdot t_1^{u_1} t_2^{u_2}\, , \,\,\,
      f_2 \,\,=\,\, (t_2 - \beta) \cdot t_1^{v_1} t_2^{v_2}\, , \,\,\,
     f_3  \,\,=\,\, (t_1 - \gamma t_2) \cdot t_1^{w_1} t_2^{w_2},$$
where $\alpha,\beta,\gamma$ are general complex numbers, and
$u_1,u_2,v_1,v_2,w_1,w_2$ are integers.
Then $\mathcal{H}$ is the arrangement of three lines in
the plane $\bR^2$ which consists of the two coordinate
axes and the main diagonal.
The images of the six one-dimensional cones of the arrangement $ 
\mathcal{H}$ under
the map $\Psi$ are spanned by the vectors
\begin{eqnarray*}
    \Psi(e_1) & = & (u_1,v_1,w_1), \\
   \Psi(e_1+e_2) & = & (u_1+u_2,v_1+v_2, w_1+w_2+1) ,\\
    \Psi(e_2) & = & (u_2, v_2, w_2) , \\
    \Psi(-e_1) & = & (-u_1-1,-v_1,-w_1-1),\\
    \Psi(-e_1-e_2) & = & (-u_1-u_2-1,-v_1-v_2-1,w_1-w_2-1), \\
    \Psi(-e_2) & = & (-u_2,-v_2-1,-w_2-1).
\end{eqnarray*}
The image of $\Psi$ consists of the six two-dimensional
cones in $\bR^3$ which are spanned by (cyclically) consecutive  
vectors in this list.
Next, the tropical surface $\cT(I)$ contains the six two-dimensional  
cones
which are spanned by the following pairs in $\bR^3$:
$$\, \{ \Psi(e_1) , e_2\}   \, , \,\,
\{  \Psi(e_1+e_2) , e_3\} \, , \,\,
\{  \Psi(e_2), e_1 \} \, , $$
$$
\{   \Psi(-e_1), e_2 \} \, , \,\,
\{   \Psi(-e_1-e_2), e_3 \} \, , \,\,
\{   \Psi(-e_2) , e_1 \}. $$
  Finally, the zero-dimensional cone $C = \{0\}$
  contributes the three two-dimensional cones which are
  spanned by the pairs  of standard basis vectors
  $$  \{ e_1, e_2 \} \,\,, \,\,\,
   \{ e_1, e_3 \} \,\,, \,\,\,
   \{ e_2, e_3 \} .$$
   These pairs determine a non-planar graph with $9$
   vertices and $15$ edges. It is isomorphic to the one depicted in  
Figure
   \ref{fig:berg36}. Thus $\cT(I)$ is a two-dimensional fan in $\bR^3$
   whose intersection with the unit sphere is an immersion of
   that graph.  The resulting embedded graph, which
   is dual to the Newton polytope of the
   implicit equation, depends on the numerical values
   of the exponents $u_1,u_2,v_1,v_2,w_1,w_2$.
   This example shows how our method can be applied to
  a family of implicitization problems where the exponents of
   the $f_i$ are not fixed integers but are unknowns.
   \qed
\end{Example}

\subsection{Surfaces}
\label{surfaces}
Suppose that $d=2$, so our given input is a
list of lattice polygons $P_1,\ldots,P_n$ in $\bR^2$.
The tropical variety $\mathcal{T}(I)$ is a two-dimensional
fan in $\bR^n$ which we represent by an embedded graph in the
$(n-1)$-sphere $\mathbb{S}^{n-1}$. This graph is constructed as follows.
We fix the node $e_i$ for each polygon $P_i$ which is not just a  
point, and we fix
the node $\Psi(w)$ for each vector $w$ that is an inner
normal to an edge of the polygon $P = P_1 + \cdots + P_n$
and such that $\Psi(w) \not= 0$. We identify
such a vector $\Psi(w)$ with
the point  $\Psi(w)/||\Psi(w)||$ on $\mathbb{S}^{n-1}$.

We now connect pairs of nodes by great circles in $\mathbb{S}^{n-1}$
according to the following rules:
a pair $\,\{\Psi(w),\Psi(w')\}\,$ gets connected if their
edges on $P$ are adjacent, a pair
$\,\{\Psi(w),e_i\} \,$ gets connected if
$w$ is the inner normal to an
edge of $P_i$, and a pair
$\,\{e_i,e_j\}\,$ gets connected if
  $P_i + P_j$ is two-dimensional.
Some pairs of great circles intersect
and thus create new nodes in $\mathbb{S}^{n-1}$.
This intersection need not be transversal, i.e., two
great circles may intersect in a smaller great circle.
The result of this construction is the tropical surface $\mathcal{T} 
(I)$,
represented by a graph in $\mathbb{S}^{n-1}$.

\begin{Example}
\label{ex:3triangles}
Consider the family of surfaces in $3$-space which is given by
$$ A_1 = \{ (1,0), (3,1), (2,2)\} ,\,
       A_2 = \{(-1,0), (0,-1), (0,0) \}, \,
       A_3 = \{ (2,1), (0,2), (1,3) \}. $$
To construct the graph $\cT(I)$ on the sphere $\mathbb{S}^2$, we first
draw the three nodes $e_1,e_2,e_3$ and the nine nodes
$\Psi(w)$ which are the images under $\Psi$ of the inner normals
of the $9$-gon $ P_1+P_2+P_3$. These nine directions are given by the
columns of the matrix
$$ \bordermatrix{
     & \, e_2 & \,e_3 & \,e_1 & \,  e_2 & \,e_3 &\, e_1 & \, e_2 &  
\,e_3 & \,e_1 \cr
     &   \phantom{-} 1 &     \phantom{-}  1 &   -1 &   -3 &   -7 &    
-4  &  -2  &     \phantom{-} 0 &     \phantom{-}  2 \, \cr
      &   -1   & -2  &   -2  &   \phantom{-}   0  &   \phantom{-}   0  
&     \phantom{-}  0  &    \phantom{-}  0  &  -1 &   -2 \, \cr
       &    \phantom{-}   2 &     \phantom{-}  4 &    \phantom{-}    
0  &  -2  &  -5  &  -4  &  -3  &  -2 &   -2 \,\cr
}
$$
We connect these twelve nodes with great circles as described
above, namely, by forming the $9$-cycle of these columns, by
connecting them to the $e_i$ as indicated, and by
forming the triangle $e_1,e_2,e_3$. This
creates an embedded graph in $\mathbb{S}^2$ which has
$14$ vertices, $27$ edges and $15$ regions. The Newton polytope
$Q$ of the implicit equation  is dual to this
graph, so it has $15$ vertices, $27$ edges and $14$ facets.
Using the  methods to be described in the next two sections, we
construct $Q$ metrically, and we find that $g(x_1,x_2,x_3)$
is a polynomial of degree $14$ having  $154$ terms. \qed
\end{Example}

\section{Multiplicities}

Let $\,K = \bC \{\!\{ \epsilon^\bR \}\!\} \,$ be the field of
Puiseux series in the unknown
$\epsilon$ which have real exponents and complex coefficients.  The field $K$ is algebraically closed, and it has
the surjective non-archimedean valuation
$$\, {\rm order} : K^* \rightarrow \bR\, , \,\,\,
\alpha \epsilon^u + \cdots \,\mapsto \, u . $$
Hence,
${\rm order}(u)$ is the exponent of the
$\epsilon$-monomial of lowest degree which appears with
non-zero coefficient in the Puiseux series $u = u(\epsilon)$. The  
order map extends
to vectors of length $n$ by coordinatewise application:
$$ {\rm order} : (K^*)^n \,\rightarrow\, \bR^n ,\,
(u_1,\ldots,u_n) \,\mapsto \,\bigl( {\rm order}(u_1), \ldots,{\rm  
order}(u_n) \bigr). $$
For any ideal $\,I \subset \bC[x_1,\ldots,x_n]\,$ we write
$V(I) \subset (K^*)^n$ for the variety defined by $I$
in the algebraic torus over the Puiseux series field $K$.
It is known (see e.g.~\cite{SS}) that the tropical variety equals
the image of the variety under the order map:
\begin{equation}
\label{kapranov}
  \mathcal{T}(I) \quad = \quad {\rm order} \bigl( V(I) \bigr).
\end{equation}
If $I$ is a prime ideal of dimension $d$, which we may assume in our implicitization problem, then  $\mathcal{T}(I)$
has the structure
of a pure $d$-dimensional polyhedral fan.
This fan structure is non-unique, and it can derived
by restricting the Gr\"obner fan of any homogenization of $I$.
See \cite{BJSST} for software and mathematical details.

Every maximal cone $\Gamma$ of the fan $\mathcal{T}(I)$ naturally  
comes with a multiplicity, which is a positive integer.  The {\em  
multiplicity} $M_\Gamma$ of a $d$-dimensional cone $\Gamma$ is the  
sum of multiplicities of all monomial-free minimal associate primes  
of the initial ideal $\text{in}_v(I)$ in $\bC[x_1, \ldots, x_n]$  
where $v \in \RR^n$ is any point in the relative interior of the cone  
$\Gamma$.  The multiplicities $M_\Gamma$
are an important piece of data which must also be
determined when computing a tropical variety $\mathcal{T}(I)$.
It is with these multiplicities that a tropical variety
satisfies the balancing condition which ensures
that tropical intersection numbers are independent of choices;
see \cite[\S 3]{Mi} and \cite[\S 6.4]{Ka}.

We now consider the prime ideal $I$  of algebraic relations among
the generic Laurent polynomials $f_1(t), \ldots,f_n(t)$.
Generators for the ideal $I$ are (still) unknown, but
in Theorem \ref{SolutionToYuProblem} we computed its
tropical variety $\mathcal{T}(I)$ in combinatorial terms.
In what follows we  similarly
compute the multiplicity $M_\Gamma$ for every maximal cone $\Gamma$
in a fan structure on $\mathcal{T}(I)$.
Recall that $\cT(I)$ is the union of cones of the form
\begin{equation}
\label{ourcones}
\Psi(C) \,+ \, \mathbb{R}_{\geq 0}^J,
\end{equation}
where $C$ is a cone of the normal fan of
  $P=P_1+\ldots+P_n$ and $\dim {\rm face}_w(P_K)\ge |K|$
for any subset $K \subseteq J$ for some (hence any)
vector $w$ in the relative interior of $C$.

The cones in  (\ref{ourcones}) do not generally form a fan, and they
need to be subdivided to give a fan structure on
the tropical variety $\cT(I)$. We choose a fine enough fan structure on
$\cT(I)$  so that any cone \eqref{ourcones}
is the union of some cones of $\cT(I)$.

We say that a pair $(C,J)$ {\em covers} a
  $d$-dimensional cone $\Gamma$ in $\cT(I)$ if
\eqref{ourcones} is $d$-dimensional and contains $\Gamma$.
There may be more than one pair $(C,J)$ which covers
a fixed cone $\Gamma$.  Suppose that
  $(C,J)$ covers $\Gamma$. Then the sublattice generated
by $\Psi(C \cap \bZ^d) + \bZ^J$ has rank $d$ in $\bZ^n$.
We define $\,{\rm index}(C,J)\,$ to be the index of that sublattice
in the maximal rank $d$ sublattice of $\bZ^n$ that contains it.

Let $F_j = {\rm face}_w(P_j)$ and
$F_J = {\rm face}_w(\sum_{j\in J} P_j) = \sum_{j\in J} F_j$
for some (hence any) vector $w$ in the relative interior of $C$.
  Then the  $|J|$-dimensional mixed volume
\begin{equation}
\label{MV2} {\rm MixedVolume}(F_j \,: \, j \in J)
\end{equation}
  is exactly the same as the one in equation (\ref{MV}).
This mixed volume is normalized with respect to the affine
lattice spanned by the $|J|$-dimensional polytope $F_J$.
Multiplying (\ref{MV2}) by index$(C,J)$ we obtain the {\em scaled  
mixed volume}
\begin{equation}
\label{MV3} {\rm MixedVolume}(F_j \,: \, j \in J) \cdot {\rm index}(C,J)
\end{equation}

The following theorem characterizes the multiplicities of our  
tropical variety.
\begin{Theorem}
\label{multFormula}
The multiplicity $M_\Gamma$ of a maximal cone $\Gamma$ in
the fan structure on $\cT(I)$ is the
sum of all scaled mixed volumes (\ref{MV3}) where the pair
$(C,J)$ covers~$\Gamma$.
\end{Theorem}

The proof of this theorem will be presented in \cite{STY},
along with the following

\begin{Remark}
Theorems \ref{SolutionToYuProblem} and \ref{multFormula}
remain valid, with the word ``cone'' replaced by the word  
``polyhedron'',
when the coefficients of the given Laurent polynomials $f_1,\ldots,f_n 
$ are no longer
complex numbers but lie in the Puiseux series field $K$.
\end{Remark}

To illustrate Theorem \ref{multFormula}, suppose that
  $d=2$ and consider the tropical surface  $\mathcal{T}(I)$.
It is represented by a graph as in Subsection \ref{surfaces}.
The edges of this graph are weighted with multiplicities
as follows. The weight of a pair $\,\{\Psi(w),\Psi(w')\}$
is the greatest common divisor of the $2 \times 2$-minors of the
$n \times 2$-matrix $\,(\Psi(w),\Psi(w'))$,
divided by the determinant of the $2 \times 2$-matrix $(w,w')$.
The weight of
a pair $\,\{\Psi(w), e_i\}$ is the normalized length of the edge $ 
{\rm face}_w(P_i)$ of $P_i$
times the greatest common divisor of all
coordinates of $\Psi(w)$ but the $i^{\rm th}$.
The weight of a pair $\,\{e_i,e_j\}\,$ is the mixed area
$$ {\rm MixedVolume}(P_i,P_j) \quad = \quad
{\rm area}(P_i+P_j) -
{\rm area}(P_i) - {\rm area}(P_j). $$
Now, when forming the embedded graph $\mathcal{T}(I) \subset
\mathbb{S}^{n-1}$, each edge $\Gamma$
may be covered by more than one of the great circles
created by these pairs, and we take $M_\Gamma$
to be the sum of their weights. With these multiplicities,
the graph $\mathcal{T}(I)$ is balanced.
Here is a concrete numerical example, taking
from the computer algebra literature.

\begin{figure}[htb]
\includegraphics[scale=0.5]{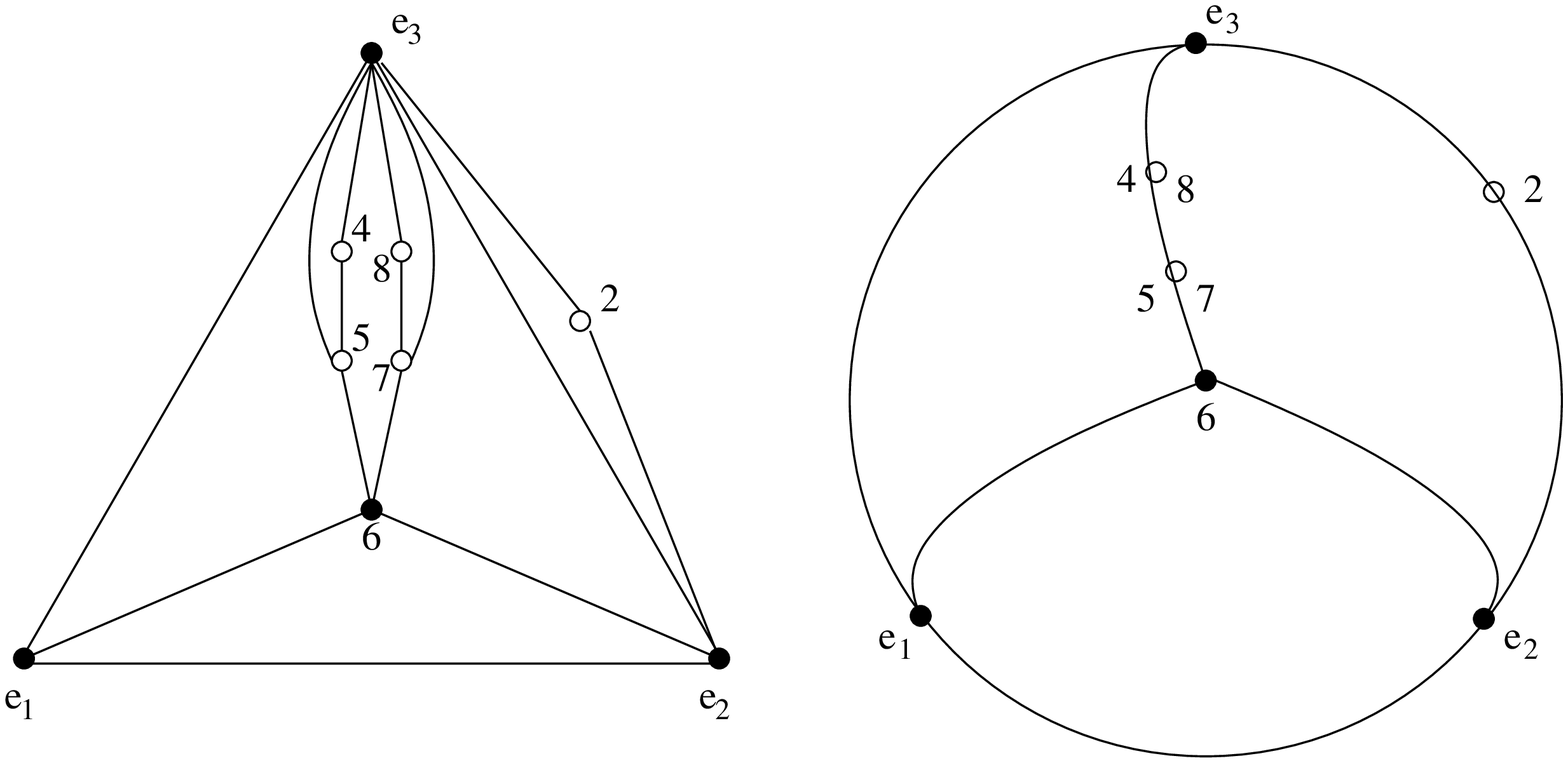}
\caption{Tropical construction of a bicubic surface in $3$-space}
\label{fig:emiris}
\end{figure}

\begin{Example}
\label{emiris}
Let $d=2$ and consider the well-known problem of implicitizing  
bicubic surfaces.
The following specific surface was discussed in \cite[Example 3.4]{EK}:
\begin{eqnarray*}
&
  f_1 (s,t) \quad  = & 3 {\bf t}^3-6 t^2+3 t+ {\bf s}^3-3 s^2+6 s-  
{\bf 1}, \\ &
  f_2 (s,t) \quad = & 3 {\bf s}^3-6 s^2+3 {\bf s}+ {\bf t}^3+ 3 {\bf  
t}, \\ &
  f_3 (s,t) \quad = & -3 {\bf s}^3 {\bf t}^3+15 s^2 t^3-15 {\bf s}  
{\bf t}^3
  -3  s^3 t^2-18  s^2 t^2 +27 s t^2   \\ &
  &       -3 {\bf t}^2+6  {\bf s}^3 {\bf t} +9  s^2 t-18  s t +3 {\bf  
t}-3 {\bf s}^2+3 {\bf s}.
\end{eqnarray*}
The vertices of the Newton polygons are indicated in boldface font.
The polygon $P_1$ is a triangle, $P_2$ is a quadrangle, and $P_3$ is  
a $7$-gon.
Their Minkowski sum $P_1+P_2+P_3$ is an octagon.
The eight primitive inner normal vectors $w = (w_1,w_2)$ of this
octagon
are listed in the first two rows of the following matrix:
$$ \bordermatrix{
    & e_1 e_2 e_3 &     e_2 e_3 & e_1 e_2 e_3  &  \,\,e_3  &  \, 
\,e_3  & e_1 e_2  & \,\,e_3  & \,\,\,e_ 2\,\, \cr
              w_1 &    1  &  1   &   0   &  -1  &  -1    &  -1    &   
\phantom{-}0   &  \phantom{-}1     \cr
              w_2 &    0  &  1   &    1  &   \phantom{-}1 &  \phantom 
{-}0   &    -1  &  -1    &  -1   \cr
        \Psi(w)_1 &     0 &   0 &   0  &  -3 &   -3 &   -3  &   -3  
&   -3 \cr
        \Psi(w)_2 &    0 &   1 &   0  &  -3 &  -3  &   -3  &  -3  &   
-3 \cr
        \Psi(w)_3 &   0 &   1  &  0  &  -2 &   -3 &   -6  &  -3  &    
-2} $$
The last three rows contain the coordinates of their
images $\Psi(w)$ under the tropical morphism $\Psi$.
We now construct the graph as in Subsection \ref{surfaces}.
  The first and third inner normal vector
is mapped to zero under $\Psi$, so they do not contribute
to the tropical surface $\mathcal{T}(I)$. The remaining
six columns $2,4,5,6,7,8$ and the unit vectors
$e_1,e_2,e_3$ form the following  $15$ great circles
on the sphere $\mathbb{S}^2$, as shown in Figure \ref{fig:emiris}:
$$
\begin{matrix}
   {\rm pair}    &     {\rm index}  &   {\rm mixed} \, {\rm volume} &  
{\rm on} \, {\rm edge} \\
\{e_1,e_2\} &              1           &   9 & \{ e_1, e_2\} \\
\{e_1, e_3\} &             1           & 18 & \{e_1, e_3 \} \\
\{e_2 ,e_3\} &             1           & 17 & \{e_2, e_3 \} \\
\{2, e_2\}      &             1           & 1   & \{e_2, e_3\} \\
\{2, e_3\}      &             1           & 1   & \{e_2, e_3 \} \\
\{6, e_1\}      &             3             & 3 & \{6, e_1 \}\\
\{6, e_2\}      &             3           & 3 &   \{6, e_2\} \\
\{4, e_3\}      &             3           & 1 &  \{6, e_3\} \\
\{8,e_3\}      &              3          & 1 & \{6, e_3\} \\
\{5, e_3 \}      &            3             & 2 & \{6, e_3\} \\
\{7,e_3\}      &              3            & 2 & \{6, e_3\} \\
\{ 4 , 5 \}      &               3         & 1 & \{6, e_3\} \\
\{  5, 6 \}      &               9         & 1 & \{6, e_3\} \\
\{  6, 7 \}      &               9         & 1& \{6, e_3\} \\
\{  7, 8 \}      &               3          & 1 & \{6, e_3\}
\end{matrix} $$
There is considerable overlap
among these great circles. For instance, the points $4$ and $8$
coincide, as do the points $5$ and $7$, and they all lie on
the great circle between $e_3$ and $6$. Likewise,
the point $2$ lies on the great circle between $e_2$ and $e_3$.
Thus the embedded graph $\mathcal{T}(I) \subset \mathbb{S}^2$ is the  
complete
graph $K_4$, and, by adding up the contributions of each
index times mixed volume, we get the corresponding weights:
$$
\begin{matrix}
\,{\rm edge} \,\,\, \Gamma  & &  \{e_1,e_2\} & \{e_1,e_3\} & \{e_2,e_3 
\} & \{6,e_1\} & \{6 ,e_2\} & \{6,e_3\} \\
{\rm weight}\,\, M_\Gamma  & & 9 & 18  & 18 & 9 & 9 & 18
\end{matrix}
$$

{F}rom this table we can now determine the Newton polytope $Q$
of the implicit equation $g(x_1,x_2,x_3)$. The general
method for this is explained in the next section.
In this example we find that the proposed Newton polytope
is the tetrahedron
$$ Q \quad = \quad \bigl\{\, (u_1,u_2,u_3) \in \bR^3 \,: \,  
u_1,u_2,u_3 \geq  0
\,\,\,\hbox{and} \,\,\, u_1 + u_2 + 2 \cdot u_3 \leq 18 \,\bigr\} .$$
The number of lattice points in $Q$ equals $715$.
What we can conclude at this point is that $Q$ would be the Newton  
polytope
of the implicit equation $g$ if $f_1,f_2,f_3$ were replaced by  
polynomials
with the same support but with generic coefficients.

However, the coefficients of the specific polynomials $f_1,f_2,f_3$
we took from \cite[Example 3.4]{EK} are not generic.
For instance, the equations $f_2(s,t) = f_3(s,t) = 0$
have only $15$ solutions $(s,t) \in (\bC^*)^2$, which is less
than the number $17$ predicted by the mixed volume.
Yet, it turns out that $Q$ is the correct
Newton polytope and all $715$ possible monomials appear
with non-zero coefficients in the implicit equation:
\begin{eqnarray*}
& \!\!\! g(x_1,x_2,x_3) \,\,= \,\,
  387420489 x_1^{18}
+387420489 x_2^{18}
-18014398509481984 x_3^9 \,+\, \cdots \cdots \\ & \,\,
-12777985432959891776936639829 x_2^2 + \cdots
-3707912273492242256259566313 .
\end{eqnarray*}
We included the $x_2^2$-term because it has the coefficient of  
largest absolute value. \qed
\end{Example}
\section{From tropical variety to Newton polytope to implicit equation}

In the previous sections we constructed the tropical variety
$\cT(I)$ along with the multiplicity $M_\Gamma$ for each
maximal cone $\Gamma$. Assuming $I = \langle g \rangle$
to be a principal ideal, we now show how this information
reveals the Newton polytope $Q$ of $g$.
The polynomial $g$ is then recovered from
$Q$ and the $f_i$'s using numerical linear algebra.
If $V(I)$ is not a hypersurface
then the role of $Q$ will be played by the
Chow polytope.

\subsection{Generic hypersurfaces}
\label{sec:genhyp}

The hypersurface case, when $n= d+1$, has recieved the most attention
in the computer algebra literature \cite{CGKW, EK, SY}. Suppose we
are given $d+1$ Laurent polynomials  $f_0,f_1,\ldots,f_d$ whose
coefficients are generic relative to their supports
$A_0,A_1,\ldots,A_d \subset \bZ^d$.  Their Newton polytopes
are denoted by $P_0,P_1,\ldots,P_d$. Note that we shifted
indices by one. The prime ideal of algebraic relations among
the $f_i$ is denoted, as before, by $\, I \subset \bC[x_0,x_1, 
\ldots,x_d]$.
There is a simple combinatorial criterion for when the ideal
$I$ is principal, i.e., when the variety $V(I)$
parametrized by $(f_0,f_1,\ldots, f_d)$ is actually a hypersurface
in $\bC^{d+1}$.

\begin{Proposition}
\label{isprincipal}
The prime ideal $I$ is principal if and only if
there exist points $a_i \in A_i$
such that the $d \times (d+1)$-matrix
$(a_0,a_1,\ldots,a_d)$ has maximal rank $d$.
\end{Proposition}

\begin{proof}
This can be seen by writing down the
$d \times (d+1)$-Jacobian matrix
$\bigl({\partial f_i}/{\partial t_j} \bigr)$,
and using the fact that the coefficients of the $f_i$ are generic.
\end{proof}

We now assume that the condition in Proposition
\ref{isprincipal}  is satisfied, and
we let $g = g(x_0,x_1,\ldots,x_d)$ denote the
unique (up to scaling) generator of the ideal $I$.
Let $Q \subset \bR^{d+1}$ be the Newton polytope
of the irreducible polynomial $g$. Then the
tropical hypersurface $\cT(I) = \cT(g)$
is the union of the codimension one cones
in the normal fan of $Q$, and the multiplicity of
a maximal cone $\Gamma$ in any fan structure on  $\cT(I)$ equals to
lattice length of the edge of $Q$ whose normal cone contains $\Gamma$.

Suppose we constructed the tropical hypersurface $\mathcal{T}(I)$
using Theorem \ref{SolutionToYuProblem}, and
we computed the multiplicities $M_\Gamma$ of each
cone $\Gamma$ in $\mathcal{T}(I)$ using Theorem \ref{multFormula}.
These data, combined with the requirement that $Q$ lies
in the non-negative orthant and intersects
each of the $d+1$ coordinate hyperplanes,
determine the polytope $Q$ uniquely.
Hence we can construct the desired Newton polytope
$\,Q \subset \bR^{d+1}\,$
combinatorially from the
given Newton polytopes $P_0,P_1,\ldots,P_d \subset \bR^d$.

The previous paragraph almost solves the problem stated in \cite{SY}.
The only shortcoming is that the description of the set $\cT(I)$
given in Theorem \ref{SolutionToYuProblem} does not come with
a nice fan structure. We saw this in
Examples \ref{ex:3triangles} and \ref{emiris}.
In our view, the following theorem provides a better solution,
as it does not require the knowledge of any fan structure.
All that is needed are the integers $M_\gamma = M_\Gamma > 0$
where $\gamma$ is a smooth point on $\cT(I)$
and $\Gamma$ is any sufficiently small
relatively open cone in $\cT(I)$ which contains $\gamma$.
But these numbers $M_\gamma$ are computed using
the formula in Theorem \ref{multFormula},
simply by replacing ``$(C,J)$ covers $\Gamma$''  with
``$(C,J)$ covers $\gamma$''.

The following theorem gives a formula for each coordinate of each vertex
of the Newton polytope $Q$. It does not rely on any particular
fan structure on $\cT(I)$.

  \begin{Theorem} \label{getinimono}
For a generic vector $v \in \bR^{d+1}$, the $i^{\rm th}$ coordinate  
of the vertex face$_v(Q)$ equals the number of intersection points,
each counted with its intersection multiplicity, of the
tropical hypersurface $\mathcal{T}(I)$ with the
halfline $\,v + \bR_{\geq 0} \,e_i$.
\end{Theorem}

This result and its generalization to larger codimension,
to be stated in Theorem \ref{AliciaEva} below, is due
to Dickenstein, Feichtner
and the first author, and it appears in the paper
\cite{DFS} on {\em Tropical Discriminants}.
We consider the conjunction of Theorems
\ref{SolutionToYuProblem}, \ref{multFormula} and
\ref{getinimono} to be a satisfactory solution
to the problem stated in  \cite{SY}.

What remains is for us to explain the meaning of the
term ``intersection multiplicity'' in
the statement of Theorem \ref{getinimono}.
   Since $v$ is generic, each intersection
of the tropical hypersurface  $\mathcal{T}(I)$ with
the halfline $\, v + \bR_{\geq 0} e_i\, $ occurs
in a smooth point $\gamma$ of
$\mathcal{T}(I)$. The {\em intersection multiplicity}
  at $\gamma$ is the product of the multiplicity $M_\gamma$
and the index of the  finite index sublattice
$\,\bZ e_i +( \bR \Gamma \cap \bZ^{d+1}) \,$ of $\bZ^{d+1}$.
  In \cite[\S 2]{DFS}, these two factors were called
the {\em intrinsic} and {\em extrinsic multiplicities} of the  
intersection.

\subsection{Computational issues}
\label{subsec:CI}
Theorems \ref{SolutionToYuProblem} and \ref{multFormula} provide
an algorithm for computing the
tropical variety $\mathcal{T}(I)$ and the multiplicities $M_\Gamma$
for each maximal cone $\Gamma$ in a fan structure of $\mathcal{T}(I)
$. The dimension of the
fan $\mathcal{T}(I)$ coincides with the dimension
of the variety $V(I)$ it represents.
When $V(I)$ is a hypersurface, characterized combinatorially
by Proposition \ref{isprincipal}, we compute a list of all
vertices of the Newton polytope $Q$ by means of
Theorem \ref{getinimono}. All of these computations
require only polyhedral geometry and linear algebra
(over the integers), but they do not require any
methods from computational commutative algebra
(Gr\"obner bases, resultants).

We are in the process of developing an implementation
of these algorithms in a software package called
{\tt TrIm} (Tropical Implicitization).
A very preliminary test implementation already exists.
It is a {\tt perl} script that incorporate the software
packages {\tt Polymake}
\cite{GJ}, {\tt BBMinkSum} \cite{Hu}, {\tt Mixed Volume Library} \cite
{EC}, {\tt Maple}, and {\tt Matlab}.  We use {\tt BBMinkSum} to
compute Minkowski sums of the input Newton polytopes $P_i$ and {\tt Polymake}
to find the face structure of  $P = P_1 + \cdots + P_n$.  Then for
each cone $C$ in the normal fan of $P$, we look for the subsets $J
\subset \{1,\dots,n\}$ with $|J| + {\rm dim}(C) = d$ such that $\Psi
(C) + \bR_{\geq 0}^J$ is in $\cT(I)$.  This amounts to doing linear
algebra using {\tt Matlab}.  We then use {\tt Mixed Volume Library}
and {\tt Maple} to compute the scaled mixed volumes \eqref{MV3} for
pairs $(C,J) \subset \cT(I)$.

The input to {\tt TrIm} is either a list of Laurent polynomials or
their Newton polytopes (in {\tt Polymake} format).  The output is a
description of $\cT(I)$ in the form of a list of pairs $(C,J)$ as
described before.  Notice that the output is not the same as the
output of {\tt Gfan} \cite{BJSST} because we do not get a fan
structure from the description in Theorem \ref{SolutionToYuProblem}.
Using this output, we can find the vertices of the Newton polytope
using Theorem \ref{getinimono} or the vertices of the Chow polytope
using Theorem \ref{AliciaEva}.

In the hypersurface case, the exponents of monomials with non-zero coefficients in the implicit equation $g$ all
lie in the Newton polytope $Q$.  After computing the vertices of $Q$,
enumerating the integer points in $Q$ gives us a list of all possible
monomials in $g$ with indeterminate coefficients.
Finally, we apply numerical linear algebra to compute the
indeterminate coefficients of the polynomial $g$.

This is done as follows.
Recall that we are given a parametrization
$(f_0, \dots, f_d)$ of the hypersurface $V(I)$.
We can thus pick any point $\tau \in (\bC^*)^d\,$ and substitute
its image  $\,(x_0,\ldots,x_d) =  (f_0(\tau), \dots, f_d(\tau))
\in V(I)\,$ into  the equation $g(x_0, \dots, x_d) = 0$.
This gives us one linear equation for the coefficients of $g$.
We now pick a different point $\tau'$ in $(\bC^*)^d\,$
to get a second equation, and so on. In this manner,
we can generate a system of linear equations whose
solution space is one-dimemensional, and is spanned
by the vector of coefficients of the desired
implicit equation $g$.

In our first experiments, we found that this
linear algebra problem seems to be numerically stable
and efficiently solvable in {\tt Matlab} if we
use vectors $\tau$ whose coordinates $\tau_i$
are unitary numbers, i.e., complex numbers of modulus $1$.
But there are other possible schemes for
generating and solving these linear systems.
For instance, the authors of
\cite{CGKW} advocate the use of integral operators
in setting up linear equations for the coefficients
on the implicit equation $g$. We believe
that implicitization is a fruitful direction of further
study in numerical linear algebra.

\subsection{Non-generic hypersurfaces}

The following proposition ensures that we can apply
{\tt Tr}opical {\tt Im}plicization even if the
given polynomials $f_i(t)$ do not have generic coefficients.
If there is a unique irreducible relation among the given
$f_i$, then our linear algebra method will find that
relation, up to a polynomial multiple.

\begin{Proposition}
\label{specialCoeffs}
Let $f_0,\ldots,f_d  \in \CC[t_1^{\pm 1},\ldots,t_d^{\pm  1}]$
be any Laurent polynomials
whose ideal of algebraic relations is principal, say
$I = \langle g \rangle$, and $P_i \subset \bR^d$ the
Newton polytope of $f_i$. Then the polytope $Q$, which is
constructed combinatorially from $P_0,\ldots,P_d $ as in
Subsection \ref{sec:genhyp}
contains a translate of the Newton polytope of $g$.
\end{Proposition}

\begin{proof}
We introduce a family of Laurent polynomials
$\,f_0^\epsilon(t), f_1^\epsilon(t), \ldots, f_d^\epsilon(t)\,$
with coefficients in $\bC(\epsilon)$ which are
generic for their support and which satisfy
$\,{\rm lim}_{\epsilon \rightarrow 0} f_i^\epsilon(t) = f_i(t) \,$
for all $i$. There exists an irreducible polynomial
$g^\epsilon(x_0,x_1,\ldots,x_d)$ with coefficients in
$\bC (\epsilon)$ which vanishes if we replace each $x_i$
by $f^\epsilon(t)$. Removing common factors
of $\epsilon$ from the terms of $g^\epsilon$, we may assume that
$$ \hat g(x) \, := \, {\rm lim}_{\epsilon \rightarrow 0} \,g^\epsilon 
(x) $$
is not the zero polynomial. The Newton polytope of $\hat g$
is contained in $Q$. Now, we take the limit
for $\epsilon \rightarrow 0$ in the identify
$$ g^\epsilon \bigl( f_0^\epsilon(t), f_1^\epsilon(t), \ldots, f_d^ 
\epsilon(t) \bigr)
\,\equiv  \, 0 \quad {\rm in} \,\, \bC(\epsilon) [t_1^{\pm 1}, 
\ldots,t_d^{\pm  1}] . $$
This implies that $\hat g(x)$ is an algebraic relation among
$f_0(t), f_1(t), \ldots, f_d(t)$, so the irreducible polynomial
$g(x)$ is a factor of $\hat g(x)$. This means that the Newton polytope
of $g$ is a Minkowski summand of the Newton polytope of $\hat g$.
This implies that a translate of the Newton polytope of $g$ lies in $Q$.
\end{proof}

\begin{Example}
Let $d = 2$ and consider the three homogeneous quadrics:
\begin{eqnarray*}
&
  f_0 (s,t) \quad  = &   \phantom{-}s^2+s t  -2 t^2  \\ &
  f_1 (s,t) \quad = &     \phantom{-}s^2- 2 st + t^2  \\ &
  f_2 (s,t) \quad = &     -2s^2+st  +t^2
  \end{eqnarray*}
  Here,  $\,g = x_0 + x_1 + x_2 $, so its Newton polytope
  is the triangle $ \Delta = {\rm conv} \{e_0,e_1,e_2\}$. The
Newton polytope of the implicit equation for
a parametric surface given by three general quadrics
is twice that triangle:  $Q = 2 \cdot \Delta$.
This shows that $\Delta$ is not a subpolytope of $Q$
but a  (non-unique) translate of $\Delta$ is a subpolytope of $Q$.
\qed
\end{Example}

The approach suggested by Proposition \ref{specialCoeffs}
is to not worry at all whether the given polynomials have
generic coefficients or special coefficients. The idea is that this
will not be relevant until the very end, when the numerical
linear algebra detects that the solution space to the
linear system of equations for the coefficients of $g(x)$
  is larger than expected.
This should work fine when the $f_i$
are not too far from the generic case. However, for
the kind of special morphisms $f$ which typically
arise in  algebraic geometry and its applications
(for instance, in statistics), this approach
is likely to fail in practice. In such cases,
one needs to be more clever in designing
suitable compactifications  $Y$ of the
very affine variety $Y^0$ of Section 2;
see~\cite{STY}.

\subsection{Lower-dimensional varieties}

We now consider the case when the parametric variety
$V(I) \subset \bC^n$ and its tropicalization
$\cT(I) \subset \bR^n$ do not have codimension one.
For technical reasons, we here assume that the $f_i$
are homogeneous polynomials of the same degree,
so that $I$ is a homogeneous prime ideal
in $\bC[x_1,\ldots,x_n]$.
Here the role of the Newton polytope $Q$
is played by the {\em Chow polytope} $\, {\rm Chow}(I)$.
This is the lattice polytope  in $\bR^n$ which is defined as follows.
For a generic vector $v \in \bR^n$,
the $i^{\rm th}$ coordinate of the vertex
$\,{\rm face}_v({\rm Chow}(I))\,$ of $\,{\rm Chow}(I)\,$
is the sum of $\mu_P({\rm in}_v(I))$  where $P$ runs over all
monomial primes that contain the variable $x_i$ and
$\mu_P({\rm in}_v(I))$  denotes the multiplicy  of
the monomial ideal ${\rm in}_v(I)$ along $P$.  The Chow polytope
was introduced and studied in \cite{KSZ}.  If $V(I)$ has codimension
one, then
${\rm Chow}(I)$ is precisely the Newton polytope
of the irreducible generator $g$ of $I$.

The following result from \cite[Theorem 2.2]{DFS} shows
that in the generic case we can construct $\,{\rm Chow}(I)\,$ from
the tropical variety $\cT(I)$ and its multiplicities
$M_\gamma$ in the same way as
constructing the Newton polytope of the implicit equation of a  
hypersurface.
Let $I$ be any prime ideal of height $c$ in $\mathbb{C}[x_1,\ldots,x_n]$
and  $v $ a generic vector in $\mathbb{R}^n$,
so that the initial ideal ${\rm in}_v(I)$ is generated by monomials.

\begin{Theorem}[DFS]
\label{AliciaEva}
   A prime  ideal  $\,P = \langle \,x_{i_1}, \ldots, x_{i_c} \rangle \,$
is a minimal prime of ${\rm in}_v(I)$ if and
   only if the tropical variety $\mathcal{T}(I)$ meets
    the cone $\,v+ \mathbb{R}_{\geq 0}\{e_{i_1},\ldots, e_{i_c}\}$.
  The number of intersections, each counted with its intersection  
multiplicity,
  coincides with  the multiplicity $\mu_P({\rm in}_v(I))$ of the  
monomial ideal ${\rm in}_v(I)$   along  $P$.
\end{Theorem}

Geometrically, the vertices of ${\rm Chow}(I)$ corresponds to the  
toric degenerations of the algebraic cycle underlying $I$.  All points in ${\rm Chow}(I)$ have the same coordinate sum,
namely, the degree of the projective variety of $I$.
Thus we get a combinatorial rule of computing the
degree of the image of any projective morphism. Moreover,
if $I$ is homogeneous with respect to some
multigrading, then  Theorem \ref{AliciaEva} gives us
a formula also for the {\em multidegree} of $I$.
The multidegree is a multivariate polynomial
which refines the classical notion of degree,
and which has received much attention recently in
algebraic combinatorics through the work of Knutson and Miller \cite 
{KM}.

\end{document}